\documentclass[11pt]{amsart}
\usepackage[margin=1in]{geometry}
\usepackage{float}
\usepackage{amssymb,amsmath,mathtools,booktabs,hyperref}
\newtheorem{theorem}{Theorem}[section]
\newtheorem{proposition}[theorem]{Proposition}
\newtheorem{lemma}[theorem]{Lemma}

\theoremstyle{definition}

\newtheorem{remark}[theorem]{Remark}
\newcommand{\F}{\mathbb{F}}
\newcommand{\Fq}{\mathbb{F}_q}
\newcommand{\Fb}{\overline{\mathbb{F}}_2}
\newcommand{\GR}{\mathrm{GR}}
\newcommand{\T}{\mathrm{T}}
\newcommand{\Tr}{\mathrm{Tr}}
\newcommand{\diff}{\delta}
\newcommand{\needlines}[1]{\par\vskip 0pt plus #1\baselineskip\penalty-200\vskip 0pt plus -#1\baselineskip}
\DeclareMathOperator{\Res}{Res}
\hypersetup{pdftitle={Explicit permutations of (Z/4Z)\string^m with uniformly bounded differential uniformity},
  pdfauthor={Hiroki Fukui}}
\DeclareMathOperator{\wt}{wt}

\title[Bounded differential uniformity on $(\mathbb Z/4\mathbb Z)^m$]{Explicit permutations of $(\mathbb Z/4\mathbb Z)^m$\\ with uniformly bounded differential uniformity}
\author{Hiroki Fukui}
\thanks{Research Institute of Criminal Psychiatry / Sex Offender Medical Center, Tokyo, Japan;
Department of Neuropsychiatry, Kyoto University, Kyoto, Japan. E-mail: \texttt{fukui@somec.org}.}
\address{Research Institute of Criminal Psychiatry / Sex Offender Medical Center, Tokyo, Japan;
Department of Neuropsychiatry, Kyoto University, Kyoto, Japan}
\email{fukui@somec.org}

\begin{document}
\begin{abstract}
Symmetric primitives over Galois rings need nonlinear components whose differentials with respect to ring addition are
controlled. Let $R=\GR(4,m)$ and $q=2^m$, so that $(R,+)\cong(\mathbb Z/4\mathbb Z)^m$ has order $N=q^2$. For every $m\ge1$ we
construct an explicit permutation $G_c$ of $R$ all of whose differentials with a nonzero input difference have probability at most
$36/N$, whereas every permutation polynomial of $R$ has such a differential of probability at least $2/\sqrt N$
(R{\o}njom--Sandrib). The bound holds for the fixed permutation itself, not on average over keys. Equivalently, $\Delta(G_c)\le36$
against $\Delta\ge2\sqrt N$ for permutation polynomials; since every permutation has $\Delta\ge2$, the least differential
uniformity of a permutation of $(\mathbb Z/4\mathbb Z)^m$ lies between $2$ and $36$ for every $m$. In Teichm\"uller digits $(a,b)\in\Fq^2$ and for $c\in\Fq$ of
absolute trace $1$,
\[
G_c(a,b)=\bigl(L,\;(L^2+L+c)\,a\bigr),\qquad L=(a^2+a+c)\,b .
\]
The map is $F_c\circ J\circ F_c$, where $F_c(a,b)=(a,(a^2+a+c)b)$ is induced by a polynomial over $R$ and $J$ exchanges the
digits: a single digit exchange between two polynomial layers yields a bound independent of $m$.
The proof reduces the difference equations, in each of three types of directions, to one-variable polynomials whose degrees do
not depend on $m$. For general maps $F\circ J\circ F$ with $F(a,b)=(p(a),h(a)+d(a)b)$, where $p,d,h$ have degrees bounded
independently of $m$, $p'\ne0$ and $d$ has no zeros, a Jacobian identity in characteristic $2$ bounds the differentials with
nonzero low input and output differences independently of $m$; the remaining directions need a separate argument, which we give
for an explicit family ($\Delta\le440$, resp.\ $510$, for odd $m$).
\end{abstract}
\maketitle

\section{Introduction}

For a map $G$ on a finite abelian group $(R,+)$ put
$\diff_G(\alpha,\beta)=\#\{x\in R: G(x+\alpha)-G(x)=\beta\}$ and $\Delta(G)=\max_{\alpha\ne0,\beta}\diff_G(\alpha,\beta)$.
Over $\F_2^n$ the group is elementary abelian, and power maps such as $x^3$ are permutations with $\Delta=2$ for every odd $n$.
Symmetric primitives over the integer rings $\mathbb Z_{2^k}$ and over Galois rings, studied for use in multiparty computation and
homomorphic encryption \cite{GMHORW,BGORS}, raise the same question for groups that are not elementary abelian. The simplest family
of this kind is $(\mathbb Z/4\mathbb Z)^m$, the additive group of $R=\GR(4,m)$, and there polynomial constructions face a
restriction: R{\o}njom and Sandrib proved that every permutation polynomial of
$\GR(p^k,m)$ has $\Delta\ge 2p^{m(k-1)}$ \cite[Theorem~3]{RS} (every polynomial function has $\Delta\ge p^{m(k-1)}$
\cite[Proposition~2]{RS}); for $\GR(4,m)$ this is $\Delta\ge 2q=2\sqrt{|R|}$, so a permutation with smaller $\Delta$ cannot be a
polynomial function.
Splitting an element into digits, applying maps to the digits and recombining is an established design strategy. It appears in the
``cut-and-sew'' S-boxes over $\mathbb Z_q$ for composite $q$ \cite[\S4.3]{GMHORW}, where for $q=p^2$ the two digits are also
exchanged; in the Kintsugi strategy over prime fields \cite[\S2]{Monolith}; and as ``digit manipulation'' over $\mathbb Z_{p^n}$
in Beyne, Grassi, {\O}ygarden, Richterov\'a and Sandrib \cite[\S9]{BGORS}. These constructions act on cyclic groups $\mathbb Z_q$
or on $\F_p$. Here we consider the groups $(\mathbb Z/4\mathbb Z)^m$, noncyclic for $m>1$, and give an explicit family whose differential
uniformity is bounded independently of $m$.
R{\o}njom and Sandrib also classify the table-defined (non-polynomial) permutations of $\GR(4,2)$ with $\Delta=2$ \cite[\S5]{RS};
no family covering all $m$ is given there. On cyclic groups $\mathbb Z_n$ several infinite families of permutations with $\Delta=2$
are known (see \cite{Tang} and the references there), but $\mathbb Z_{4^m}$ and $(\mathbb Z/4\mathbb Z)^m$ are different groups for
$m>1$.
Permutation polynomials are subject to the lower bound $\Delta\ge2\sqrt{|R|}$ \cite{RS}, and the one-sided digit-exchange maps
$J\circ F$ of \cite{FukuiDE} attain $\Delta=q=\sqrt{|R|}$, which is optimal within that class (Theorem~\ref{thm:oneswap}).
The family constructed here has $\Delta\le36$ for every $m$, so the optimal order is constant.

\needlines{8}
\begin{theorem}\label{thm:main}
Let $m\ge1$, $q=2^m$, and $c\in\Fq$ with $\Tr_{\Fq/\F_2}(c)=1$. In Teichm\"uller digits, the map
$G_c=F_c\circ J\circ F_c$ with $F_c(a,b)=(a,D(a)b)$, $D(x)=x^2+x+c$, $J(a,b)=(b,a)$, is a permutation of $\GR(4,m)$ and
\[
\Delta(G_c)\le 36 .
\]
\end{theorem}

Let $D_m$ be the least $\Delta$ over all permutations of $(\mathbb Z/4\mathbb Z)^m$. Theorem~\ref{thm:main} and
Proposition~\ref{prop:two} give $2\le D_m\le36$ for every $m$, so $D_m=\Theta(1)$, whereas every permutation polynomial has
$\Delta\ge2\sqrt{|R|}$. In terms of probabilities, with $N=|R|=q^2$, every differential of $G_c$ with a nonzero input difference has probability at most
$36/N$, which is of the optimal order $N^{-1}$. For example, for $m=8$ the ring $R$ has $65{,}536$ elements (a $16$-bit state):
every permutation polynomial of $R$ has a differential of probability at least $512/65{,}536=1/128$, while every differential of $G_c$
with a nonzero input difference has probability at most $36/65{,}536$, more than $14$ times smaller.
In the instances of Table~\ref{tab:data} ($m\le10$) the exact values of $\Delta(G_c)$ are at most $20$. Since $F_c$ is itself
induced by a polynomial (Section~\ref{ss:maps}), the polynomial lower bound is overcome by a single digit exchange between two
polynomial layers. The proof (Section~\ref{sec:proof}) is elementary, and the degrees it uses do not depend on $m$.

\medskip
\noindent\textbf{Comparison.} With $F_{p,d,h}(a,b)=(p(a),h(a)+d(a)b)$:
\begin{center}
\begin{tabular}{p{0.44\textwidth}ll}
\toprule
map & digit system & $\Delta$\\
\midrule
permutation polynomial & --- & $\ge 2q$ \cite{RS}\\
$J\circ F$ (output exchange only) & any & $\ge q$ (Thm.~\ref{thm:oneswap})\\
$F_c\circ J\circ F_c$ & Teichm\"uller & $\le 36$ for all $m$ (Thm.~\ref{thm:main})\\
$F\circ J\circ F$, $p=x^3$, $d=x^4+x^2+1$,\newline $h\in\{0,x^2,x^3,x^5\}$, $m$ odd & Teichm\"uller & $\le 440$, resp.\ $510$ (Thm.~\ref{thm:family})\\
literal exchange $J$\newline ($=F\circ J\circ F$ with $p=x$, $d=1$, $h=0$) & Teichm\"uller & $=q$ (Remark~\ref{rem:scope})\\
literal exchange $J$ & coordinate & $=q^2/2$ (Prop.~\ref{prop:Jsigma})\\
$F\circ J\circ F$, $1$-periodic $h,d$ & coordinate & $\ge q/2$ (Prop.~\ref{prop:sigma})\\
\bottomrule
\end{tabular}
\end{center}
The class $J\circ F$ in the second line contains the maps studied in \cite{FukuiDE} (coordinate digits, $p$ APN, $d$ affine), for
which $\Delta=q$ is attained; Theorem~\ref{thm:oneswap} shows that $q$ cannot be improved within this class. Unlike \cite{FukuiDE}, which also gives key-averaged and multi-round
statements, we consider a single fixed permutation and all its differentials. For the families of Theorems~\ref{thm:main}
and~\ref{thm:family}, one more polynomial layer after the exchange removes the dependence on $m$.

\medskip
\noindent\textbf{Scope.} We do not determine $D_m$, and we do not consider other attacks (linear, integral, algebraic) or
implementation cost. Open questions are collected in Section~\ref{sec:open}.

\section{Setting}

\subsection{Digits}
$R=\GR(4,m)=W_2(\Fq)$. Every $x\in R$ is uniquely $x=\T(a)+2\T(b)$ with $a,b\in\Fq$, where $\T$ is the Teichm\"uller lift; we write
$x=(a,b)$. The identity $\T(a)+\T(u)=\T(a+u)+2\T(\sqrt{au})$ gives
\begin{equation}\label{eq:add}
(a,b)+(u,v)=(a+u,\;b+v+\sqrt{au}),\qquad -(a,b)=(a,\;a+b),
\end{equation}
hence
\begin{equation}\label{eq:sub}
(y_0,y_1)-(z_0,z_1)=(y_0+z_0,\;y_1+z_0+z_1+\sqrt{y_0z_0}).
\end{equation}
The \emph{coordinate} digit system fixes a basis $e_0=1,e_1,\dots,e_{m-1}$ of $\Fq$ over $\F_2$ whose first vector is $1$ (for
definiteness, the polynomial basis $1,\omega,\dots,\omega^{m-1}$ of $\Fq=\F_2(\omega)$) and lifts coordinates to $\{0,1\}\subset\mathbb Z/4$.
The resulting map depends on the basis; all statements below about this system are for a basis with $e_0=1$. Then
\eqref{eq:add} holds with the carry $\sqrt{au}$ replaced by the coordinatewise product $a\odot u$.

We apply the same formulas in both digit systems; they give different permutations of $R$. The change of digits has the form $C(a,b)=(a,b+\gamma(a))$, under which $F_{p,d,h}$ becomes
$F_{p,d,h'}$ with $h'(a)=h(a)+d(a)\gamma(a)+\gamma(p(a))$, and $CJC^{-1}$ is in general not the literal exchange $J$.

\subsection{The maps}\label{ss:maps}
For $p,d,h\in\Fq[x]$ with $p$ a permutation of $\Fq$ and $d$ without zeros in $\Fq$, $F=F_{p,d,h}$ is a permutation of $\Fq^2$, hence
of $R$; so are $J\circ F$ and $F\circ J\circ F$. Theorem~\ref{thm:main} is the case $p=x$, $h=0$, $d=x^2+x+c$.

We call $F_{p,d,h}$ a polynomial layer because it is induced by a polynomial over $R$. For $\mathsf F\in R[X]$ with reduction
$\bar{\mathsf F}\in\Fq[X]$, the identity $\mathsf F(y+2z)=\mathsf F(y)+2z\,\mathsf F'(y)$ shows that $\mathsf F$ acts on digits as
$(a,b)\mapsto(\bar{\mathsf F}(a),\,h_0(a)+\bar{\mathsf F}'(a)\,b)$, where $h_0(a)$ is the high digit of $\mathsf F$ at $(a,0)$; conversely,
every triple of functions $p,d,h$ arises in this way, by interpolation with the factor $X^q-X$ and a correction $2\hat H$ of the high
digit. This is \cite[Proposition ``Coordinate form and realizability'']{FukuiDE}, stated there for coordinate digits; the argument
uses only that $2w\in R$ depends on $w\bmod 2$ alone, and holds verbatim for Teichm\"uller digits. The exchange $J$ is not of this form. Since every permutation polynomial of $R$ has $\Delta\ge2q$ \cite{RS}, Theorem~\ref{thm:main}
shows that for $m\ge5$ the permutation $G_c$ is not induced by any polynomial over $R$.

\section{\texorpdfstring{Output exchange alone: $\Delta\ge q$}{Output exchange alone}}

\begin{theorem}\label{thm:oneswap}
Let $p,h:\Fq\to\Fq$ be arbitrary and $d:\Fq\to\Fq^\times$. In either digit system, $\Delta(J\circ F_{p,d,h})\ge q$.
\end{theorem}
\begin{proof}
Teichm\"uller digits: $B=J\circ F$ is $B(a,b)=(h(a)+d(a)b,\;p(a))$. Since $|\Fq^\times|<q$, $d$ is not injective; choose $a_0$ and
$u\ne0$ with $d(a_0+u)=d(a_0)$, and $v$ with $h(a_0+u)+h(a_0)+d(a_0)(v+\sqrt{a_0u})=0$. For $x=(a_0,b)$ and $\alpha=(u,v)$ the low
digits of $B(x+\alpha)$ and $B(x)$ coincide, so by \eqref{eq:sub}
$B(x+\alpha)-B(x)=(0,\,p(a_0+u)+p(a_0))$ for all $q$ values of $b$. The coordinate case is identical with $a_0\odot u$ in place of
$\sqrt{a_0u}$ (for the maps of \cite{FukuiDE}, $\Delta(B)=q$ exactly).
\end{proof}

\begin{proposition}\label{prop:two}
Every map $G:(\mathbb Z/4\mathbb Z)^m\to(\mathbb Z/4\mathbb Z)^m$, $m\ge1$, has $\Delta(G)\ge2$. Consequently, for every $m$,
$2\le\min_G\Delta(G)\le36$, the minimum over all permutations.
\end{proposition}
\begin{proof}
Let $\alpha\ne0$ with $2\alpha=0$ and put $D(x)=G(x+\alpha)-G(x)$. Then $D(x+\alpha)=G(x)-G(x+\alpha)=-D(x)$. If some $x$ has
$D(x)=\beta$ with $2\beta=0$, then $D(x+\alpha)=\beta$ as well, so $\diff_G(\alpha,\beta)\ge2$. Otherwise $D$ misses the $2^m\ge2$
elements of order at most $2$, so $D$ is not surjective and some value is taken twice.
\end{proof}
For permutations the bound $\Delta\ge2$ is immediate in any finite group of order $>2$, since $x\mapsto G(x+\alpha)-G(x)$ never takes the
value $0$; Proposition~\ref{prop:two} extends it to arbitrary maps.

\section{\texorpdfstring{Proof of Theorem~\ref{thm:main}}{Proof of the main theorem}}\label{sec:proof}

\emph{Permutation.} The image of $x\mapsto x^2+x$ on $\Fq$ is the trace-zero hyperplane, so $\Tr(c)=1$ gives $D(a)\ne0$ for all $a$;
$F_c$ multiplies $b$ by a nonzero scalar on each fibre and is bijective. Explicitly
\[
G_c(a,b)=\bigl(L,\;D(L)\,a\bigr),\qquad L=D(a)\,b .
\]

\emph{Difference equations.} Fix $\alpha=(u,v)\ne0$ and $\beta=(s,t)$, and write $a=A^2$, $u=U^2$ (Frobenius is bijective on $\Fq$,
so $A$ runs over $\Fq$ and $\sqrt{au}=AU$). By \eqref{eq:add}, $x+\alpha=(a+u,\,b+v+AU)$, so
\[
L=D(a)\,b,\qquad L'=D(a+u)\,(b+v+AU)
\]
are the low digits of $G_c(x)$ and $G_c(x+\alpha)$. By \eqref{eq:sub}, $G_c(x+\alpha)-G_c(x)=\beta$ is equivalent to
\begin{equation}\label{eq:low}
L'+L=s,
\end{equation}
\begin{equation}\label{eq:high}
t=D(L')(a+u)+L+D(L)\,a+\sqrt{LL'} .
\end{equation}
Given \eqref{eq:low}, $D(L')=D(L)+s^2+s$ and $\sqrt{LL'}=\sqrt{L^2+sL}=L+\sqrt{sL}$, so \eqref{eq:high} becomes
$t=u\,D(L)+(s^2+s)(a+u)+\sqrt{sL}$, and squaring (Frobenius is injective) gives the equivalent condition
\begin{equation}\label{eq:one}
u^2\bigl(L^4+L^2+c^2\bigr)+(s^4+s^2)(a^2+u^2)+sL+t^2=0 .
\end{equation}
We count pairs $(A,b)\in\Fq^2$ satisfying \eqref{eq:low} and \eqref{eq:one}; the map $(A,b)\mapsto x=(A^2,b)$ is a bijection.

\emph{Case A: $u\notin\{0,1\}$.} Put $\kappa=u^2+u\ne0$, so $D(a+u)=D(a)+\kappa$. Then \eqref{eq:low} reads
$\kappa b+D(a+u)(v+AU)+s=0$, which determines $b$ from $A$, and
\[
L=N(A)/\kappa,\qquad N(A)=D(a)\,s+D(a)D(a+u)(v+AU),
\]
a polynomial in $A$ of degree $9$ with leading coefficient $U$. Substituting into \eqref{eq:one} and multiplying by $\kappa^4$ gives
\[
u^2\bigl(N^4+\kappa^2N^2+c^2\kappa^4\bigr)+(s^4+s^2)(A^4+u^2)\kappa^4+sN\kappa^3+t^2\kappa^4=0,
\]
a polynomial in $A$ of degree $36$ with leading coefficient $u^2U^4=U^8\ne0$ (all other terms have degree at most $18$). Hence at most
$36$ values of $A$, each with one $b$: $\diff\le36$. This includes $s=0$.

\emph{Case B: $u=1$.} Then $U=1$ and $D(a+1)=D(a)$, so \eqref{eq:low} reads $D(A^2)(v+A)+s=0$, a monic polynomial of degree $5$ in
$A$: at most $5$ values. For each, \eqref{eq:one} is a monic polynomial of degree $4$ in $L$, and $b=L/D(a)$. Hence $\diff\le20$.

\emph{Case C: $u=0$, $v\ne0$.} Then \eqref{eq:low} reads $v\,D(a)=s$; there is no solution unless $s\ne0$, and then $D(a)=s/v$ has at most
$2$ roots $a$. For each, \eqref{eq:one} becomes $(s^4+s^2)a^2+sL+t^2=0$, linear in $L$ with coefficient $s\ne0$. Hence $\diff\le2$.

The three cases cover all $\alpha\ne0$, so $\Delta(G_c)\le36$. \qed

\begin{remark}
The argument holds for all $m$ and uses neither an APN property nor intersection theory. The constant comes from Case~A: the
degree $36=4\cdot9$ of $u^2N^4$. In the instances of Table~\ref{tab:data} the observed values are at most $20$ and are not monotone
in $m$ ($10$ for $m=5$, $9$ for $m=6$).
\end{remark}

\begin{remark}[Evaluation]
On Teichm\"uller digits, $G_c$ costs two general multiplications, two squarings and four additions in $\Fq$; given the output
$(L,y)$, the inverse $a=y/D(L)$, $b=L/D(a)$ costs two inversions, two multiplications and two squarings. Passing between $x\in R$ and its digits
requires $a=x\bmod 2$ and $\T(a)$, i.e.\ $m$ squarings in $R$ (since $\T(a)=\tilde a^{\,q}$ for any lift $\tilde a$). In settings
where elements of $R$ are held as shares modulo $4$, the reduction $x\bmod2$ is a bit extraction whose cost depends on the protocol;
we do not analyse it here.
\end{remark}

\section{\texorpdfstring{A general method for $F\circ J\circ F$}{A general method for FJF}}

For general $(p,d,h)$ the low-digit equation is no longer linear in $b$, and a uniform bound needs an argument about common
components. Let $G=F\circ J\circ F$, $\alpha=(u,v)\ne0$, $\beta=(s,t)$, $a=A^2$, $u=U^2$, and
\[
L=h(a)+d(a)b,\quad L'=h(a+u)+d(a+u)(b+v+AU),
\]
\[
Z_0=p(L),\ Z_1=h(L)+d(L)p(a),\quad Y_0=p(L'),\ Y_1=h(L')+d(L')p(a+u).
\]
By \eqref{eq:sub}, after squaring the high digit,
\begin{equation}\label{eq:E}
\diff_G(\alpha,\beta)=\#\{(A,b)\in\Fq^2:\ E_1=E_2=0\},\quad
E_1=Y_0+Z_0+s,\quad E_2=(Y_1+Z_0+Z_1+t)^2+Y_0Z_0 .
\end{equation}

\begin{lemma}\label{lem:jac}
Put $\Psi=Y_0+Z_0$ and $\Phi=(Y_1+Z_0+Z_1)^2+Y_0Z_0$. Then
$\det\partial(\Psi,\Phi)/\partial(A,b)=\Psi\cdot p'(L)\,p'(L')\cdot d(a)\,d(a+u)\cdot U$.
\end{lemma}
\begin{proof}
In characteristic $2$ squares have zero derivative, so $\Phi$ may be replaced by $Y_0Z_0$. The pair $(\Psi,Y_0Z_0)$ consists of the
elementary symmetric functions of $(Z_0,Y_0)$, with Jacobian $Z_0+Y_0=\Psi$; next $\partial(Z_0,Y_0)/\partial(L,L')=p'(L)p'(L')$.
Since $L$ involves $A$ only through $a=A^2$, $\partial L/\partial A=0$, $\partial L/\partial b=d(a)$, $\partial L'/\partial b=d(a+u)$,
$\partial L'/\partial A=d(a+u)U$, so $\partial(L,L')/\partial(A,b)=d(a)d(a+u)U$.
\end{proof}

\begin{theorem}\label{thm:general}
Let $p,d,h\in\Fq[x]$ with $p'\ne0$ and $d$ without zeros in $\Fq$, and let $K=\Fb$. If $u\ne0$ and $s\ne0$, every common component of
$E_1$ and $E_2$ in $K^2$ is a line $A=c_0$ with $c_0\notin\Fq$, and $\diff_G(\alpha,\beta)\le\deg E_1\cdot\deg E_2$.
\end{theorem}
\begin{proof}
Let $\Gamma$ be an irreducible common component. On $\Gamma$, $\Psi\equiv s$ and $\Phi\equiv t^2$, so both differentials vanish on the
tangent line at smooth points and the Jacobian of Lemma~\ref{lem:jac} vanishes on $\Gamma$. As $\Psi=s\ne0$ and $U\ne0$, $\Gamma$ lies in
$\{p'(L)=0\}$, $\{p'(L')=0\}$, $\{d(a)=0\}$ or $\{d(a+u)=0\}$. The last two are unions of lines $A=c_0$ with $c_0^2$ or $c_0^2+u$ a root of
$d$, so $c_0\notin\Fq$.

Suppose $p'(L)=0$ on $\Gamma$, so $L\equiv r$ for a root $r$ of $p'$. If $\Gamma$ is vertical, $A\equiv c_0$ while $b$ varies along $\Gamma$, so
$L=h(c_0^2)+d(c_0^2)b\equiv r$ forces $d(c_0^2)=0$ and $c_0\notin\Fq$. If $\Gamma$ is not vertical, $A$ is transcendental over $K$ on $\Gamma$ and
$b=(r-h(A^2))/d(A^2)$, so the function field is $K(\Gamma)=K(A)$ and
\[
L'|_\Gamma=E(A^2)+U\,d(A^2+u)\,A,\qquad E\in K(A^2).
\]
Since $K(A)=K(A^2)\oplus A\,K(A^2)$ and $U\,d(A^2+u)\ne0$, $L'|_\Gamma\notin K$; hence $\Psi|_\Gamma=p(r)+p(L'|_\Gamma)$ is not constant,
contradicting $\Psi\equiv s$. The case $p'(L')=0$ is symmetric: then $L|_\Gamma=E'(A^2)+U\,d(A^2)\,A$ with $E'\in K(A^2)$.

Thus $C=\gcd(E_1,E_2)$ is a product of vertical lines without $\Fq$-points, every $\Fq$-solution is an isolated point of
$V(E_1/C)\cap V(E_2/C)$, and affine B\'ezout gives the bound.
\end{proof}

The second-stage high digit enters $E_2$ only inside the square, so Theorem~\ref{thm:general} holds for every $h$. The bound
depends on $m$ only through the degrees of $p,d,h$, and is uniform in $m$ when these degrees are bounded.

\begin{remark}\label{rem:scope}
Theorem~\ref{thm:general} bounds $\diff_G(\alpha,\beta)$ only for $u\ne0$, $s\ne0$, and assumes $p'\ne0$. The directions with $u=0$ or
$s=0$ need a separate argument, and without one the differential uniformity of $F\circ J\circ F$ need not be bounded independently of
$m$. For example, $p=x$, $d=1$, $h=0$ satisfy the hypotheses, but $F=\mathrm{id}$ and $F\circ J\circ F=J$. For $\alpha=(1,0)$ and
$x=(0,b)$ we have $x+\alpha=(1,b)$, and by \eqref{eq:sub}
\[
J(x+\alpha)-J(x)=(b,1)-(b,0)=(0,1)\qquad\text{for all }b\in\Fq,
\]
so $\Delta(J)\ge q$. In fact $\Delta(J)=q$: for $\alpha=(u,v)$ put $s=v+\sqrt{au}$; then $J(x+\alpha)-J(x)=(s,\,u+\sqrt{bs})$.
If $u\ne0$, the value $s$ determines $a$; then $t$ determines $b$ when $s\ne0$, and when $s=0$ we have $t=u$ with $b$ free. If $u=0$,
then $s=v\ne0$, $t$ determines $b$, and $a$ is free. The same row $\alpha=(1,0)$ gives $\Delta\ge q$ for the literal exchange in the
coordinate system, but its exact value there is $q^2/2$ (Proposition~\ref{prop:Jsigma}); the two literal exchanges are different
permutations of $R$. In
Theorem~\ref{thm:main} and Theorem~\ref{thm:family} the remaining directions are treated by hand.
\end{remark}

\begin{theorem}\label{thm:family}
Let $m$ be odd, $p=x^3$, $d=x^4+x^2+1$, $h\in\{0,x^2,x^3,x^5\}$, and $G=F\circ J\circ F$ in Teichm\"uller digits. Then
$\Delta(G)\le 440$ for $h\in\{0,x^2,x^3\}$ and $\Delta(G)\le510$ for $h=x^5$.
\end{theorem}
\begin{proof}
For $m$ odd, $p$ is bijective and $d=(x^2+x+1)^2$ has no zeros in $\Fq$.

\emph{$u\ne0$, $s\ne0$.} By Theorem~\ref{thm:general} the $\Fq$-solutions are isolated roots in $K^2$ of a system whose supports are
contained in the generic supports of $(E_1,E_2)$; the total degrees are at most $(27,84)$, resp.\ $(29,98)$, which already gives
$2268$, resp.\ $2842$. These supports contain the origin and a pure power of each variable, so the pair is cornered, and it is nice
(e.g.\ $(A^{27}+1,b^8+1)$, resp.\ $(A^{29}+1,b^{10}+1)$, has finitely many common roots). Rojas' Affine Point Theorem~II, valid over
any algebraically closed field, bounds the number of isolated affine roots, with multiplicity, of every system with these supports by
the mixed volume of the Newton polygons \cite[Affine Point Theorem~II]{Rojas}; this applies to each specialised system. The Newton polygons are
\begin{align*}
&P_1=\mathrm{conv}\{(0,0),(27,0),(25,2),(16,3),(0,3)\},\\ &P_2=\mathrm{conv}\{(0,0),(84,0),(72,8),(0,8)\}\quad(\mathrm{MV}=440),\\
&P_1=\mathrm{conv}\{(0,0),(29,0),(25,2),(16,3),(0,3)\},\\ &P_2=\mathrm{conv}\{(0,0),(98,0),(82,8),(64,10),(0,10)\}\quad(\mathrm{MV}=510),
\end{align*}
the second pair for $h=x^5$.

\emph{$u\ne0$, $s=0$.} On $\Fq$-points $p(L)=p(L')$ forces $L=L'$, and then $h$ cancels:
$E_2\iff t=u\,d(L)(a^2+ua+u^2)$. If $u\notin\{0,1\}$, put $c_1=u^4+u^2\ne0$; then $L=L'$ gives $L=N/c_1$ with
$N=d(a)h(a+u)+d(a+u)h(a)+d(a)d(a+u)(v+AU)$, of degree $17$ in $A$ with leading coefficient $U$, and the remaining equation is
non-constant of degree $72$ in $A$; $b$ is then determined. If $u=1$, then $d(a+1)=d(a)$, $L'=L$ reads $h(a+1)+h(a)+d(a)(v+A)=0$, whose
odd part $d(a)A$ is nonzero, so at most $9$ values of $A$; for each, $d(L)=t/(a^2+a+1)$ and $L\mapsto L^4+L^2$ is two-to-one. Hence
$\diff\le72$.

\emph{$u=0$, $v\ne0$.} Put $w=v\,d(a)\ne0$; in the coordinates $(a,L)$ (a bijection of $\Fq^2$),
$E_1=wL^2+w^2L+w^3+s$. If $s=0$ then $L^3=(L+w)^3$, impossible as cubing is bijective and $w\ne0$; so $s\ne0$. Reducing $E_2$ modulo $E_1$
over $K(a)$ leaves a polynomial of degree $1$ in $L$ with coefficient $s^2/w\ne0$, so each $a$ gives at most one $L$. Eliminating $L$,
$\Res_L(E_1,E_2)=d(a)^{2}g(a)$ (resp.\ $d(a)^4g(a)$ for $h=x^5$) with $g$ of degree $88$ (resp.\ $96$), nonzero leading coefficient for
$v\ne0$, and only even powers of $a$. As $\Fq$ is perfect, $g$ is the square of a polynomial of degree $44$ (resp.\ $48$), so
$\diff\le44$ (resp.\ $48$).
\end{proof}

\begin{remark}
For the doubly exchanged maps $J\circ F\circ J\circ F$ of the same family, an analogous Jacobian identity holds,
$\det\partial(Y_1+Z_1,Y_1Z_1)/\partial(A,b)=(Y_1+Z_1)\,g(L)g(L')\,d(a)d(a+u)U$ with $g=h'+d'\cdot p$; whether these maps also have $\Delta$
bounded independently of $m$ is left open (Section~\ref{sec:open}).
\end{remark}

\section{The coordinate digit system}\label{sec:coord}

\begin{proposition}\label{prop:sigma}
In the coordinate digit system with $e_0=1$, if $d(a+1)=d(a)$ and $h(a+1)=h(a)$ for all $a$ (for instance $h=0$, $d=x^4+x^2+1$, or $h=0$,
$d=x^2+x+c$), then $\Delta(F\circ J\circ F)\ge q/2$.
\end{proposition}
\begin{proof}
Take $\alpha=(1,0)$. Since $1=e_0$, the carry is $a\odot1=a_0$, the first coordinate of $a$. For the
$q/2$ values of $a$ with $a_0=0$ and all $b$, $x+\alpha=(a+1,b)$ and $L'=L$, so $G(x+\alpha)-G(x)=\bigl(0,\;d(L)(p(a+1)+p(a))\bigr)$.
These $q^2/2$ inputs land in at most $q$ values of $\beta$.
\end{proof}

For $h=0$ and $d=x^4+x^2+1$ (with $p=x^3$) the row $\alpha=(1,0)$ alone reaches $6,22,74$ for $m=3,5,7$. In general the
kernel of $a\mapsto a\odot u$ has $2^{m-\wt(u)}$ elements; in the Teichm\"uller system, by contrast,
$a\mapsto\sqrt{au}$ is injective for $u\ne0$. In the coordinate system the carry is a linearised polynomial of degree up to $2^{m-1}$, so the bounded-degree
arguments above do not apply; whether $F\circ J\circ F$ has bounded $\Delta$ there for some $h$ that is not $1$-periodic is open.

\begin{proposition}\label{prop:Jsigma}
In the coordinate digit system with $e_0=1$, the literal exchange $J_\sigma(a,b)=(b,a)$ has $\Delta(J_\sigma)=q^2/2$.
\end{proposition}
\begin{proof}
For $\alpha=(u,v)$ put $s=v+a\odot u$. Then $J_\sigma(x+\alpha)-J_\sigma(x)=(s,\,u+b\odot s)$. For a target $(s,t)$, the set of $a$
with $a\odot u=s+v$ has $0$ or $2^{m-\wt(u)}$ elements, and the set of $b$ with $b\odot s=t+u$ has $0$ or $2^{m-\wt(s)}$ elements;
hence $\diff(\alpha,(s,t))\le2^{2m-\wt(u)-\wt(s)}$. If $u\ne0$ then $\wt(u)\ge1$; if $u=0$ then $s=v\ne0$ and $\wt(s)\ge1$. In both
cases $\wt(u)+\wt(s)\ge1$, so $\diff\le q^2/2$. For
$\alpha=(0,e_0)$ and $\beta=(e_0,0)$, every $a$ and every $b$ with $b_0=0$ is a solution, so the bound is attained.
\end{proof}
For $m=1$ both literal exchanges have $\Delta=2$; for $m>1$ they differ.

\section{Data}
\begin{table}[H]
\caption{Exact $\Delta(G_c)$ for the tested parameters. Ring model: $\GR(4,m)=\mathbb Z_4[X]/(f)$ with ring addition only and digits
recovered by table lookup, for all $c$ with $\Tr(c)=1$ when $m\le5$ and for two values of $c$ when $m=6$. Digit formulas
\eqref{eq:add}--\eqref{eq:sub} (all $N-1$ directions, with a coverage check): one value of $c$ for each $m\le10$ and two for $m=8$;
these agree with the ring model for $m\le6$. The field polynomial and the values of $c$ are listed in the supplementary material.
For each $m$, the value was the same for every tested $c$.}\label{tab:data}
\small
\begin{tabular}{rrrrrrrrrrr}
\toprule
$m$ & 1 & 2 & 3 & 4 & 5 & 6 & 7 & 8 & 9 & 10\\
$q$ & 2 & 4 & 8 & 16 & 32 & 64 & 128 & 256 & 512 & 1024\\
\midrule
$\Delta(G_c)$ & 2 & 4 & 4 & 6 & 10 & 9 & 12 & 14 & 16 & 20\\
\bottomrule
\end{tabular}
\end{table}

\section{Open problems}\label{sec:open}
\begin{enumerate}
\item Determine $D_m$. By Table~\ref{tab:data} and \cite[\S5]{RS}, $D_1=D_2=2$. Are there permutations of
$(\mathbb Z/4\mathbb Z)^m$ with $\Delta=2$ for every $m$?
\item How far can the constant $36$ be lowered, for $G_c$ or for other explicit families? For the tested parameters with
$m\le10$ (Table~\ref{tab:data}) the values of $\Delta(G_c)$ are at most $20$.
\item For families $(p,d,h)$ of bounded degree, Theorem~\ref{thm:general} bounds the differentials with $u\ne0$, $s\ne0$
independently of $m$. For which such families are the remaining directions also bounded independently of $m$? Does the same hold for the doubly exchanged maps
$J\circ F\circ J\circ F$?
\item Can $\Delta(F\circ J\circ F)$ remain bounded independently of $m$ in the coordinate digit system for some $h$ that is not
$1$-periodic (Section~\ref{sec:coord})?
\end{enumerate}

\appendix
\section{Computational checks}\label{app:cert}
The proofs above are complete without computation. The following scripts (supplementary material; Magma, and one C program with its runner) check the stated
identities and degrees with \texttt{assert}, and recompute the table.
\begin{itemize}
\item \texttt{t4\_v2\_bounded36\_symbolic.m}: Theorem~\ref{thm:main}, from the unreduced equations \eqref{eq:low}--\eqref{eq:high}:
in Case A, $N$ of degree $9$ with leading coefficient $U$ and the eliminated equation of degree $36$ with leading coefficient
$U^8/\kappa^4$; in Case B, the low equation monic of degree $5$ in $A$ and the high equation monic of degree $4$ in $L$; in Case C,
the low equation $vD(a)+s$ and the high equation linear in $L$ with coefficient $s$.
\item \texttt{t5\_ringmodel\_v2.m} (parameters \texttt{mlo}, \texttt{mhi}, \texttt{cmax}): $\Delta(G_c)$ in the ring model, using ring
addition only.
\item \texttt{gc\_ddt\_v2.c} with \texttt{run\_gc\_v2.sh}: $\Delta(G_c)$ from the digit formulas, in parallel shards; the runner accepts a
run only if every shard exits with status $0$ and the recorded number and sum of directions of each shard equal those of its residue
class $\alpha\equiv k\pmod K$.
\item \texttt{c1\_form1\_assert.m}: Lemma~\ref{lem:jac} as a polynomial identity (the family of Theorem~\ref{thm:family} and random
$(p,d,h)$), and the mixed volumes $440$ and $510$.
\item \texttt{t3\_u0\_form1\_v2.m}: the $u=0$ elimination in Theorem~\ref{thm:family} (exponent of the removed factor $d(a)$, degree and
leading coefficient of $g$, even exponents, and the linear remainder).
\item \texttt{t1\_J\_counterexample.m}: $\Delta(J)=q$ in the ring model for $m\le5$ (Teichm\"uller digits, Remark~\ref{rem:scope});
\texttt{r05\_t1\_digit\_swap.m}: both literal exchanges for $m\le5$, with $\Delta=q$ (Teichm\"uller) and $\Delta=q^2/2$ (coordinate,
Proposition~\ref{prop:Jsigma}).
\item \texttt{t6\_sigma\_row.m}: the row $\alpha=(1,0)$ in the coordinate system (Proposition~\ref{prop:sigma}).
\end{itemize}
The supplementary material contains the scripts, their raw outputs, a manifest with SHA-256 hashes, and instructions for
reproducing each run. It is archived with the paper at \url{https://doi.org/10.5281/zenodo.23875215} and at\\
\url{https://github.com/hirokifukui/z4m-bounded-differential-uniformity}.

\section*{Funding}
This research received no external funding.

\section*{Use of AI tools}
The author used ChatGPT (OpenAI) and Claude (Anthropic) as research assistants, for exploring constructions and proof strategies,
checking arguments, writing verification scripts, and drafting and editing the text. The author reviewed and revised all content
and takes full responsibility for it. The numerical results were obtained with the supplied C and Magma programs.

\end{document}